# On the matrix form of second-order linear difference equations

M.I. Ayzatsky[1]

National Science Center Kharkov Institute of Physics and Technology (NSC KIPT),
610108, Kharkov, Ukraine

Results of research of possibility of transformation of a difference equation into a system of the first-order difference equation are presented. In contrast to the method used previously, an unknown grid function is split into two new auxiliary functions, which have definite properties. Several examples show that proposed approach can be useful in solving different physical problems.

## 1 Introduction

"Yet other and more general types of sets are needed in order to determine the most natural context for defining difference equations and improve their applicability to solving problems in biology, economics, computer science and other fields." ([1], preface)

It is common knowledge that a difference equation of order $k$ may be transformed in a standard way to a system of $k$ first-order difference equations. For example, a second-order[2] difference equation

$$y_{k+1} + a_k y_k + b_k y_{k-1} = f_k. \tag{1.1}$$

can be rewritten as

$$Y_k = T_k Y_{k-1} + F_k, \tag{1.2}$$

where

$$Y_k = \begin{pmatrix} x_k \\ y_k \end{pmatrix}, \ x_k = y_{k+1}, \ T_k = \begin{pmatrix} -a_k & -b_k \\ 1 & 0 \end{pmatrix}, \ F_k = \begin{pmatrix} f_k \\ 0 \end{pmatrix}. \tag{1.3}$$

This is a pure mathematic transformation that does not take into account the possible physical meaning of the solution of the equation (1.1). For example, the process of wave propagation can be described by the equation (1.1) which can be a grid approximation of the wave equation [2,3,4,5,6,7,8,9,10,11,12] or represent relations of the amplitudes of the field expansion [13,14,15,16,17,18,19,20,21]. It is desirable that the components of the vector $Y_k$ represent some physical notions. For the case of wave propagation very useful notions are "forward and backward waves" that form the general field. The problem of creating a special field distribution is also needed in such approach [15,18,22,23] Moreover, following question can be formulated. Is the transformation of a $k$-order linear difference equation to a system of $k$ first-order linear difference equations unique or there are several ones?

In this paper we present results of our research of possible transformation of the equation (1.1) into the equation (1.2).

I would like to note that before writing this note I overviewed great number of mathematical papers and books (some of them are [1,24,25,26,27,28,29,30,31,32,33,34,35,36,37,38]), and did not find similar result. However, there are so many works in the field of difference equations that I could easily miss the important contribution. I apologize in advance, if I did not refer on the works that have relevance to the topic under discussion.

---

[1] M.I. Aizatskyi, N.I.Aizatsky; aizatsky@kipt.kharkov.ua

[2] It is a well-known fact that second-order differential problems are very often encountered in the applications, especially among those derived from physics. A more frequent appearance of second-order problems is also true in difference equations.

## 2 Transformation a second-order linear difference equation

We seek a solution[3] of difference equation (1.1) as

$$y_k = y_k^{(1)} + y_k^{(2)}.$$ (1.4)

where $y_k^{(1)}$, $y_k^{(2)}$ are the new unknown grid functions.

By introducing two unknowns $y_k^{(1)}$, $y_k^{(2)}$ instead of the one $y_k$, we can impose an additional condition. Let us assume that

$$y_{k+1} = \rho_k^{(1)} y_k^{(1)} + \rho_k^{(2)} y_k^{(2)},$$ (1.5)

where $\rho_k^{(1)}$ and $\rho_k^{(2)}$ ($\rho_k^{(1)} \neq \rho_k^{(2)}$) are the given numbers. If we define $y_{k+1}$ and $y_k$, then from (1.4) and (1.5) we can find unique solutions for $y_k^{(1)}$ and $y_k^{(2)}$

$$y_k^{(1)} = \frac{y_{k+1} - \rho_k^{(2)} y_k}{\left(\rho_k^{(1)} - \rho_k^{(2)}\right)}$$
$$y_k^{(2)} = -\frac{y_{k+1} - \rho_k^{(1)} y_k}{\left(\rho_k^{(1)} - \rho_k^{(2)}\right)}.$$ (1.6)

It proves the correctness of our representation (1.4) - (1.5).

We would like to emphasize that the sequences $\rho_k^{(1)}$ and $\rho_k^{(2)}$ are the arbitrary ones, and we do not impose a condition that the new grid functions $y_k^{(1)}$, $y_k^{(2)}$ are the solutions of equation(1.1).

For Cauchy problems it is necessary to give two initial values for $y_1$, $y_2$ and consider the equation (1.1) and the expression (1.5) with $k \geq 3$, the expression (1.4) with $k \geq 2$. The initial conditions for $y_k^{(1)}$, $y_k^{(2)}$ are:

$$y_2^{(1)} = \frac{f_2 - y_1 b_2 - \left(\rho_2^{(2)} + a_2\right) y_2}{\left(\rho_2^{(1)} - \rho_2^{(2)}\right)}$$
$$y_2^{(2)} = -\frac{f_2 - y_1 b_2 - \left(\rho_2^{(1)} + a_2\right) y_2}{\left(\rho_2^{(1)} - \rho_2^{(2)}\right)}.$$ (1.7)

Representations (1.4) - (1.5) have no analogy in the theory of differential equations and smooth functions. Indeed, in the case when $y_k$, $y_{k+1}$ and $y_{k+2}$ have close values ($y_{k+2} \approx y_{k+1} \approx y_k$), $y_k^{(1)}$ and $y_{k+1}^{(1)}$ ($y_k^{(2)}$ and $y_{k+1}^{(2)}$) may have very different values

$$y_{k+1}^{(1)} - y_k^{(1)} \approx y_k \left[ \frac{1 - \rho_{k+1}^{(2)}}{\left(\rho_{k+1}^{(1)} - \rho_{k+1}^{(2)}\right)} - \frac{1 - \rho_k^{(2)}}{\left(\rho_k^{(1)} - \rho_k^{(2)}\right)} \right].$$ (1.8)

From(1.1), (1.4) and (1.5) it follows that we have system:

$$y_{k+1}^{(1)} + y_{k+1}^{(2)} = \rho_k^{(1)} y_k^{(1)} + \rho_k^{(2)} y_k^{(2)},$$
$$\left(\rho_{k+1}^{(1)} + a_{k+1}\right) y_{k+1}^{(1)} + \left(\rho_{k+1}^{(2)} + a_{k+1}\right) y_{k+1}^{(2)} = f_{k+1} - b_{k+1}\left(y_k^{(1)} + y_k^{(2)}\right).$$ (1.9)

This system can be easily transformed into the normal system of difference equations:

$$y_{k+1}^{(1)} = T_{k,11} y_k^{(1)} + T_{k,12} y_k^{(2)} + F_{k+1}$$
$$y_{k+1}^{(2)} = T_{k,21} y_k^{(1)} + T_{k,22} y_k^{(2)} - F_{k+1},$$ (1.10)

or in matrix form:

---

[3] We do not make assumption that we are working with real numbers

$$\begin{pmatrix} y_{k+1}^{(1)} \\ y_{k+1}^{(2)} \end{pmatrix} = T_k \begin{pmatrix} y_k^{(1)} \\ y_k^{(2)} \end{pmatrix} + \begin{pmatrix} F_{k+1} \\ -F_{k+1} \end{pmatrix}. \tag{1.11}$$

Transmission matrix $T_k$ has following components

$$T_{k,11} = -\frac{\left[b_{k+1} + \left(\rho_{k+1}^{(2)} + a_{k+1}\right)\rho_k^{(1)}\right]}{\left(\rho_{k+1}^{(1)} - \rho_{k+1}^{(2)}\right)}. \tag{1.12}$$

$$T_{k,12} = -\frac{\left[b_{k+1} + \left(\rho_{k+1}^{(2)} + a_{k+1}\right)\rho_k^{(2)}\right]}{\left(\rho_{k+1}^{(1)} - \rho_{k+1}^{(2)}\right)}. \tag{1.13}$$

$$T_{k,21} = \frac{\left[b_{k+1} + \left(\rho_{k+1}^{(1)} + a_{k+1}\right)\rho_k^{(1)}\right]}{\left(\rho_{k+1}^{(1)} - \rho_{k+1}^{(2)}\right)}. \tag{1.14}$$

$$T_{k,22} = \frac{\left[b_{k+1} + \left(\rho_{k+1}^{(1)} + a_{k+1}\right)\rho_k^{(2)}\right]}{\left(\rho_{k+1}^{(1)} - \rho_{k+1}^{(2)}\right)}. \tag{1.15}$$

$$F_{k+1} = \frac{f_{k+1}}{\left(\rho_{k+1}^{(1)} - \rho_{k+1}^{(2)}\right)}. \tag{1.16}$$

In this paper we will consider the homogeneous difference equations. It is means that $f_k = 0$.

Now we can describe several properties of the normal system of difference equations(1.11).

From (1.12)-(1.15) it follows that we can choose the sequences $\rho_k^{(1)}$ and $\rho_k^{(2)}$ in a such way that matrix $T_k$ will be triangular or even diagonal one. The first case is realized by setting $T_{k,12} = 0$

$$T_{k,12} = -\frac{\left[b_{k+1} + \left(\rho_{k+1}^{(2)} + a_{k+1}\right)\rho_k^{(2)}\right]}{\left(\rho_{k+1}^{(1)} - \rho_{k+1}^{(2)}\right)} = 0, \tag{1.17}$$

This condition gives the non-linear second-order rational difference equation (Riccaty type[4] [30,31,35])

$$\rho_{k+1}^{(2)} = -\frac{b_{k+1}}{\rho_k^{(2)}} - a_{k+1} = -\frac{b_{k+1} + \rho_k^{(2)}a_{k+1}}{\rho_k^{(2)}}, \tag{1.18}$$

The system (1.10) takes the form

$$\begin{aligned} y_{k+1}^{(1)} &= T_{k,11} y_k^{(1)} \\ y_{k+1}^{(2)} &= T_{k,21} y_k^{(1)} + T_{k,22} y_k^{(2)} \end{aligned}, \tag{1.19}$$

From (1.19) it follows that the grid function $y_k^{(1)}$ is independent from the grid function $y_k^{(2)}$ and proportional to the product of the factors $T_{k,11}$. The equation (1.18) define the sequences $\rho_k^{(2)}$.

If $\rho_k^{(1)}$ is also the solution of the equation(1.18), but $\rho_k^{(1)} \neq \rho_k^{(2)}$, matrix $T_k$ is a diagonal and the system (1.10) takes the form:

$$\begin{aligned} y_{k+1}^{(1)} &= \rho_k^{(1)} y_k^{(1)} \\ y_{k+1}^{(2)} &= \rho_k^{(2)} y_k^{(2)} \end{aligned}. \tag{1.20}$$

---

[4] In [30], p.177 there are very useful relations between the Ricatti equation and the linear second-order difference equation

In this case $y_k^{(1)}$, $y_k^{(2)}$ are the linearly independent solutions [38] of equation(1.1). Indeed, the determinant $\Delta_2(y_k^{(1)}, y_k^{(2)})$

$$\Delta_2(y_k^{(1)}, y_k^{(2)}) = \begin{vmatrix} y_k^{(1)} & y_{k+1}^{(1)} \\ y_k^{(2)} & y_{k+1}^{(2)} \end{vmatrix} = \begin{vmatrix} y_k^{(1)} & \rho_k^{(1)} y_k^{(1)} \\ y_k^{(2)} & \rho_k^{(2)} y_k^{(2)} \end{vmatrix} = y_k^{(1)} y_k^{(2)} \left( \rho_k^{(2)} - \rho_k^{(1)} \right) \neq 0 \qquad (1.21)$$

is not zero.

If we know $y_k^{(1)}$, $y_k^{(2)}$ with some initial conditions, we can construct the solution of a boundary problem as

$$y_k = C_1 y_k^{(1)} + C_2 y_k^{(2)}, \qquad (1.22)$$

where $C_1, C_2$ are the solution of the system of linear equations (for the boundary problem of the first kind)

$$\begin{aligned} C_1 y_1^{(1)} + C_2 y_1^{(2)} = y_1 \\ C_1 y_N^{(1)} + C_2 y_N^{(2)} = y_N \end{aligned}. \qquad (1.23)$$

For a difference equation with constant coefficients $a_k = a$, $b_k = b$, the nonlinear difference equation (1.18) has two stationary points $\rho_1, \rho_2$ that are the solutions of the characteristic square equation

$$\rho^2 + a\rho + b = 0. \qquad (1.24)$$

Setting $\rho_k^{(1)} = \rho_1$, $\rho_k^{(2)} = \rho_2$ we can write

$$\begin{aligned} y_{k+1}^{(1)} = \rho_1^{k-2} y_2^{(1)}, \quad k \geq 3 \\ y_{k+1}^{(2)} = \rho_2^{k-2} y_2^{(2)} \end{aligned}. \qquad (1.25)$$

This is a well-known result [24-38].

We must note that the stationary point of the equation (1.18) can be unstable.

Triangular or diagonal systems of difference equations can be useful in many applications. One of them is finding conditions when the grid functions $y_k^{(1)}$, $y_k^{(2)}$ have the given properties.

Let's, for example, find the condition when $b_k$ and $a_k$ are not constants, but we want $\rho_k^{(1)}$ in (1.20) to be a constant $\rho_k^{(1)} = \rho_*$. From (1.18) we obtain that in this case $b_{k+1}$ and $a_{k+1}$ have to be linearly proportional

$$b_{k+1} = -\rho_* a_{k+1} - \rho_*^2, \qquad (1.26)$$

and $\rho_k^{(2)} \neq const$ have to be found from equation

$$\rho_{k+1}^{(2)} = -\frac{b_{k+1}}{\rho_k^{(2)}} - a_{k+1} = -\frac{\rho_*^2 + \rho_* a_{k+1}}{\rho_k^{(2)}} - a_{k+1}, \qquad (1.27)$$

In some cases, it is useful[5] to work with S-matrix (see, for example, [10,11,39])

$$\begin{pmatrix} y_k^{(2)} \\ y_{k+1}^{(1)} \end{pmatrix} = S_k \begin{pmatrix} y_k^{(1)} \\ y_{k+1}^{(2)} \end{pmatrix} + \begin{pmatrix} F_{k+1}^{(S)} \\ -\rho_k^{(2)} F_{k+1}^{(S)} \end{pmatrix}. \qquad (1.28)$$

Components of S-matrix for a second-order difference equation (1.1) are

$$S_{k,11} = -\frac{\left[ b_{k+1} + \rho_k^{(1)} \left( \rho_{k+1}^{(1)} + a_{k+1} \right) \right]}{\left[ b_{k+1} + \rho_k^{(2)} \left( \rho_{k+1}^{(1)} + a_{k+1} \right) \right]}. \qquad (1.29)$$

$$S_{k,12} = \frac{\left( \rho_{k+1}^{(1)} - \rho_{k+1}^{(2)} \right)}{\left[ b_{k+1} + \rho_k^{(2)} \left( \rho_{k+1}^{(1)} + a_{k+1} \right) \right]}. \qquad (1.30)$$

---

[5] It is known that S-matrix is useful when the solutions have exponential growth.

$$S_{k,21} = \frac{\left(\rho_k^{(1)} - \rho_k^{(2)}\right)b_{k+1}}{\left[b_{k+1} + \rho_k^{(2)}\left(\rho_{k+1}^{(1)} + a_{k+1}\right)\right]}. \tag{1.31}$$

$$S_{k,22} = -\frac{\left[b_{k+1} + \rho_k^{(2)}\left(\rho_{k+1}^{(2)} + a_{k+1}\right)\right]}{\left[b_{k+1} + \rho_k^{(2)}\left(\rho_{k+1}^{(1)} + a_{k+1}\right)\right]}. \tag{1.32}$$

$$F_{k+1}^{(S)} = \frac{f_{k+1}}{\left[b_{k+1} + \rho_k^{(2)}\left(\rho_{k+1}^{(1)} + a_{k+1}\right)\right]}. \tag{1.33}$$

## 2 Some physical examples

### 2.1 Wave propagation

We will consider only a few problems from the great number of ones (see [2-10] and literature cited there) to illustrate the usefulness of the considered above transformation.

We shall deal with difference equation ("diagonal tight-binding description" [11])

$$y_{k+1} + y_{k-1} - \left(2 - h^2 \varepsilon_k\right) y_k = 0. \tag{1.34}$$

that is a grid approximation of such differential equation

$$\frac{d^2 y}{d\xi^2} + \varepsilon(\xi) y = 0. \tag{1.35}$$

This differential equation is usually considered in the framework of model, in which instead of the continuous variation of the permittivity $\varepsilon(\xi)$, the piecewise constant law is introduced [4,7,8,10,11]. Matrix form that was proposed in the previous section will be similar to this approach [4,7,8,10,11] if $\rho_k^{(1)}$, $\rho_k^{(2)}$ are the solutions of the "local" characteristic equation:

$$\rho_k^2 - \left(2 - h^2 \varepsilon_k\right)\rho_k + 1 = 0. \tag{1.36}$$

$$\rho_k^{(1)} = \left(1 - \frac{h^2 \varepsilon_k}{2}\right) + \sqrt{\left(1 - \frac{h^2 \varepsilon_k}{2}\right)^2 - 1}$$
$$\rho_k^{(2)} = \left(1 - \frac{h^2 \varepsilon_k}{2}\right) - \sqrt{\left(1 - \frac{h^2 \varepsilon_k}{2}\right)^2 - 1} \tag{1.37}$$

For $\varepsilon_k = 1$

$$\rho_{vac}^{(1)} = \left(1 - \frac{h^2}{2}\right) + \sqrt{\left(1 - \frac{h^2}{2}\right)^2 - 1}$$
$$\rho_{vac}^{(2)} = \left(1 - \frac{h^2}{2}\right) - \sqrt{\left(1 - \frac{h^2}{2}\right)^2 - 1} \tag{1.38}$$

Let's consider, as an example, wave diffraction on the inhomogeneous slab in the case of permittivity $\varepsilon_k$ changing according to following law:

$$\varepsilon_k = \begin{cases} 1, \ k = 1,...,N_1 - 1 \\ 1 - 2\dfrac{k - N_1 + 0.5}{N_2 - N_1}, \ k = N_1,...,N_2 \\ 1, \ k > N_2 \end{cases}. \tag{1.39}$$

For $k > \left(N_2 + N_1\right)/2$, permittivity $\varepsilon_k$ becomes negative and electromagnetic wave does not propagate in this part of slab.

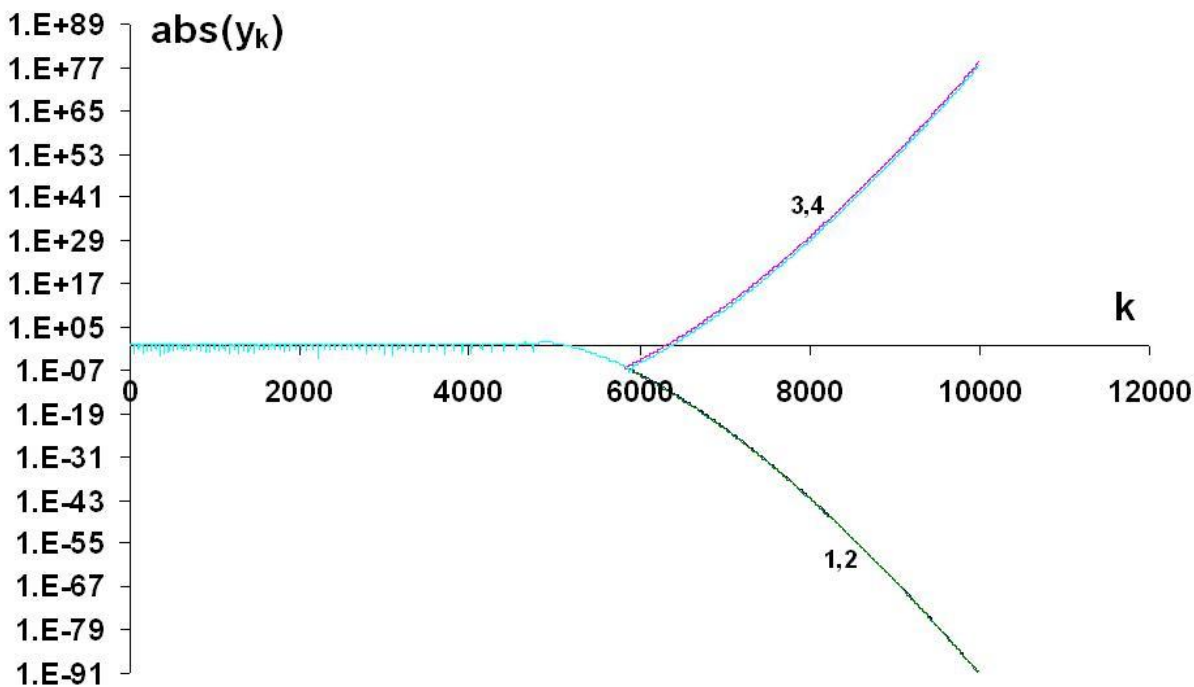

Fig. 1 Graphs of grid functions that are the solutions of the systems of the first-order difference equations (1.28) (S-matrix approach, graphs 1,2) and (1.10) (T-matrix approach, graphs 3,4). Values of the grid functions 1 and 3 were calculated for $\rho_k^{(1)}$, $\rho_k^{(2)}$ that are the solutions of the "local" characteristic equation (1.36), Values of the grid functions 2 and 4 were calculated for $\rho_k^{(1)} = const = 0.9$, $\rho_k^{(2)} = const = 1/\rho_k^{(1)}$. ($h = 2\pi/100$)

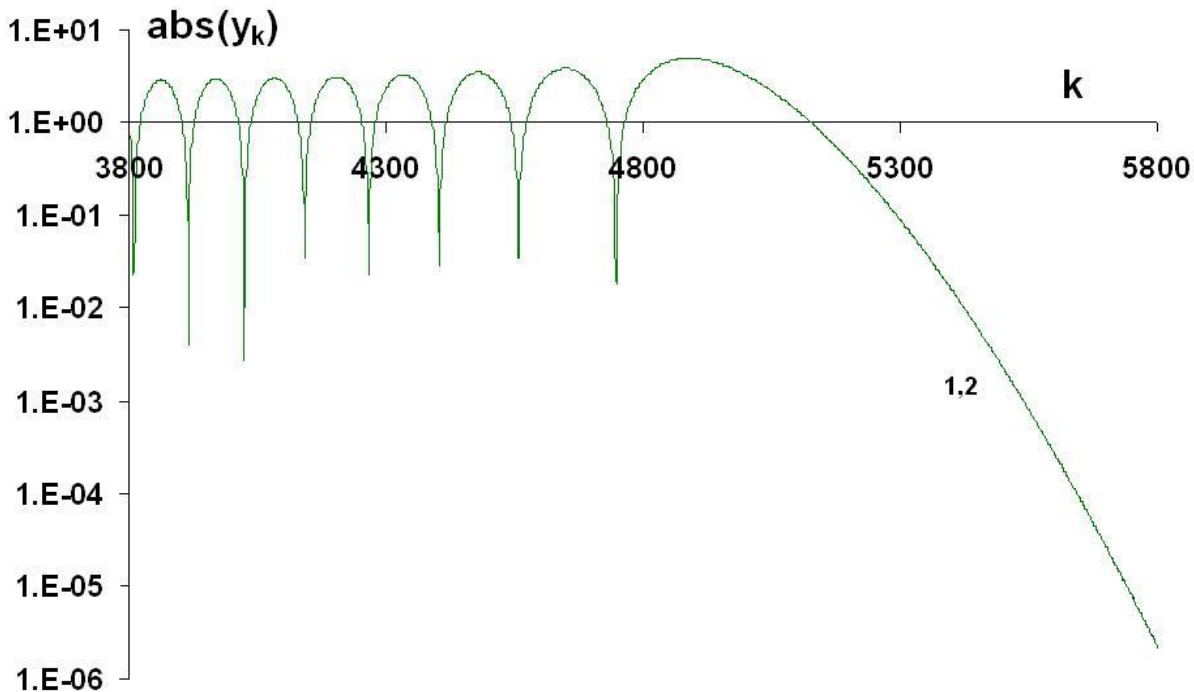

Fig. 2 Graphs of grid functions that are the solutions of the systems of the first-order difference equations (1.28) (S-matrix approach, graphs 1,2). Values of the grid function 1 were calculated for $\rho_k^{(1)}$, $\rho_k^{(2)}$ that are the solutions of the "local" characteristic equation(1.36), values of the grid functions 2 were calculated for $\rho_k^{(1)} = const = 0.9$, $\rho_k^{(2)} = const = 1/\rho_k^{(1)}$.

In Fig. 1 and Fig. 2 results of numerical solving of systems of difference equations (1.28) and (1.10), that correspond to the difference equation (1.34), are presented ( $y_k = y_k^{(1)} + y_k^{(2)}$ ). We can see that the choice of the values of $\rho_k^{(1)}$, $\rho_k^{(2)}$ does not effect on the numerical solutions of the difference equation(1.34). Moreover, as would be expected, using the transfer matrix in this problem is incorrect as it leads to divergent solutions.

The most interesting question in the using of the proposed matrix form is the properties of the partial solutions $y_k^{(1)}$, $y_k^{(2)}$ in the case when matrix $T_k$ is a diagonal.

Let's consider the simplest example of the wave diffraction on the homogeneous dielectric slab

$$\varepsilon = \begin{cases} 1, & 0 \le \xi \le \xi_1 \\ \varepsilon_2, & \xi_1 < \xi \le \xi_2 \\ 1, & \xi > \xi_2 \end{cases}. \tag{1.40}$$

Using the standard "mode matching technique" [40], we obtain analytical expression for the reflection and transmission coefficients

$$R = \frac{(\varepsilon_2 - 1) 2i \sin\left\{\sqrt{\varepsilon_2}\left(\xi_2 - \xi_1\right)\right\} \exp(i 2 \xi_1 + i \xi_2)}{\left[\left(\sqrt{\varepsilon_2} + 1\right)^2 \exp\left\{-i\sqrt{\varepsilon_2}\left(\xi_2 - \xi_1\right)\right\} - \left(\sqrt{\varepsilon_2} - 1\right)^2 \exp\left\{i\sqrt{\varepsilon_2}\left(\xi_2 - \xi_1\right)\right\}\right]}$$

$$T = \frac{4\sqrt{\varepsilon_2} \exp(i\xi_1 - i\xi_2)}{\left[\left(\sqrt{\varepsilon_2} + 1\right)^2 \exp\left\{-i\sqrt{\varepsilon_2}\left(\xi_2 - \xi_1\right)\right\} - \left(\sqrt{\varepsilon_2} - 1\right)^2 \exp\left\{i\sqrt{\varepsilon_2}\left(\xi_2 - \xi_1\right)\right\}\right]}. \tag{1.41}$$

Consider the case when the slab permittivity $\varepsilon_2 = 3 + i \times 0.03$ and $\xi_1 = 2\pi$, $\xi_2 = 11 \times 2\pi$. For this parameters $R = $-0.3207-i×6.5787E-002 ( $|R| = 0.3273$ ), $T = $-0.2185+i×0.4836 ( $|T| = 0.5306$ ). Grid approximation for the slab permittivity ( $h = 2\pi / 100$ ) is:

$$\varepsilon_k = \begin{cases} 1, k = 1, ..., N_1 & N_1 = 100 \\ 3 + i\,0.03, k = N_1 + 1, ..., N_2 & N_2 = 1100 \\ 1, \ k = N_2 + 1, ..., N_3 & N_3 = 1200 \end{cases}. \tag{1.42}$$

Using the S-matix formalism (1.28) with $\rho_k^{(1)}$, $\rho_k^{(2)}$ that are the solutions of the "local" characteristic equation(1.36), we can calculate the reflection and transmission coefficients and the values of the grid function $y_k = Y_k = y_k^{(1)} + y_k^{(2)}$ ( for $k < N_1$ and $k > N_2$ $y_k^{(1)}$ is a forward wave and $y_k^{(2)}$ is a backward wave). For $y_1^{(1)} = 1$, $y_1^{(1)} = R$ and $y_k^{(2)} = 0$, $k > N_2$ ( $N_1 = 100$ ), we obtained $R^{(Y)} = S_{11}^{\Sigma} = $-0.318-i×5.1929E-002 ( $|R^{(Y)}| = 0.3222$ ), $T^{(Y)} = S_{21}^{\Sigma} = $-0.2145+i×0.4859 ( $|T^{(Y)}| = 0.5312$ ). Comparison $R^{(Y)}$ and $T^{(Y)}$ with the exact values $R$ and $T$ shows good agreement. Graph of the grid function $Y_k$ is presented in Fig. 3

If $\rho_k^{(1)}$, $\rho_k^{(2)}$ are the solutions of the Riccaty equations (1.18), $y_k^{(1)}$ and $y_k^{(2)}$ are the linearly independent solutions. For finding these solutions we have to set the initial conditions $\rho_1^{(1)}$, $\rho_1^{(1)}$. It is convenient to choose $\rho_1^{(1)}$, $\rho_1^{(1)}$ that are the solutions of the "local" characteristic equation (1.36) with $\varepsilon_1 = 1$: $\rho_1^{(1)} = \rho_{vac}^{(1)}$, $\rho_1^{(2)} = \rho_{vac}^{(2)}$. Results of calculation the solutions of the Riccaty equations (1.18) $y_k^{(1)}$ and $y_k^{(2)}$ with the following initial conditions $y_1^{(1)} = 1$, $y_1^{(2)} = 1$ are presented in Fig. 4 and Fig. 5.

Are these functions $y_k^{(1)}$, $y_k^{(2)}$ associated with the function $Y_k$? Comparison of these grid functions (see the graphs, which are presented in Fig. 6), shows that the following relation is fulfilled

$$Y_k = y_k^{(1)} + R^{(Y)} y_k^{(2)}. \qquad (1.43)$$

We obtained the result that is expected for mathematics, but not the usual for physics.

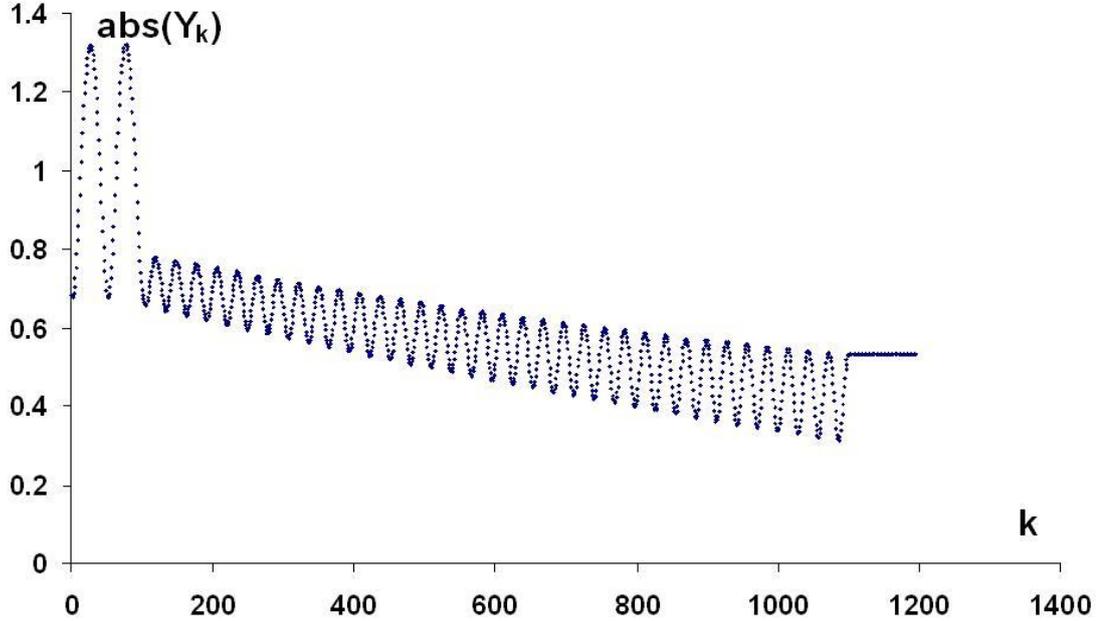

Fig. 3 Graphs of the modulus of grid function $Y_k = y_k^{(1)} + y_k^{(2)}$, where $y_k^{(1)}, y_k^{(2)}$ are the solutions of the systems of the first-order difference equations (1.28) (S-matrix approach). Matrix $T_k$ is not a diagonal.

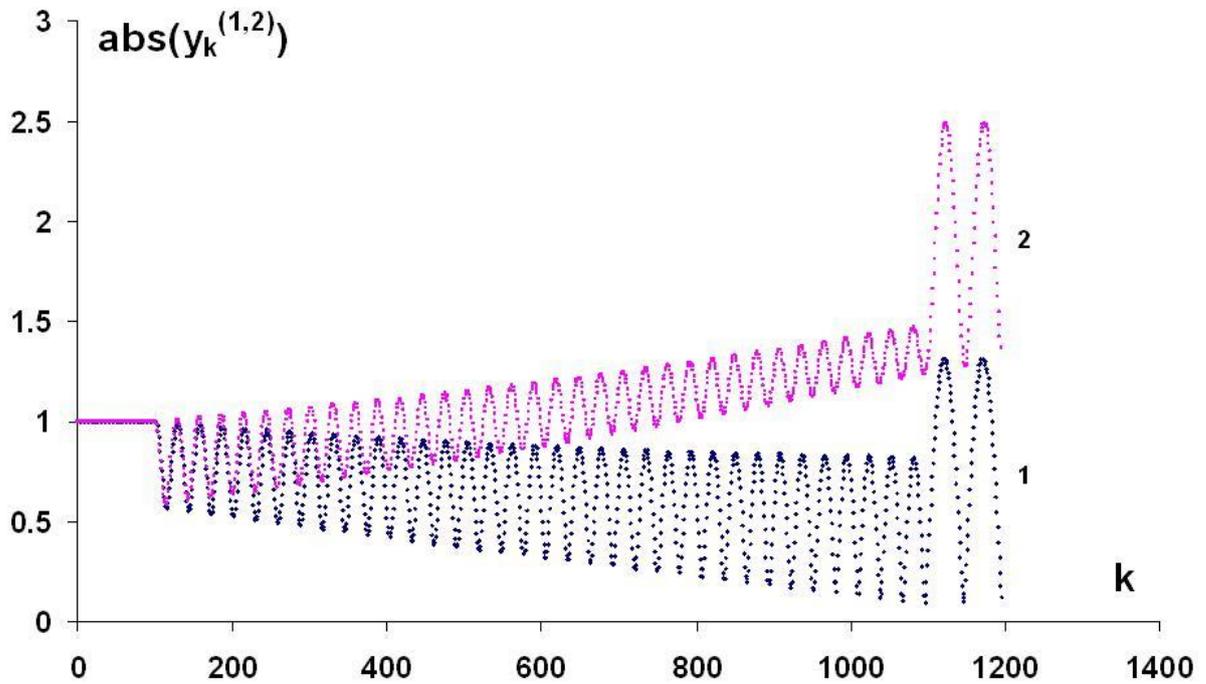

Fig. 4 Graphs of modulus of the grid functions $y_k^{(1)}, y_k^{(2)}$ for the wave diffraction on the homogeneous dielectric slab. Matrix $T_k$ is a diagonal.

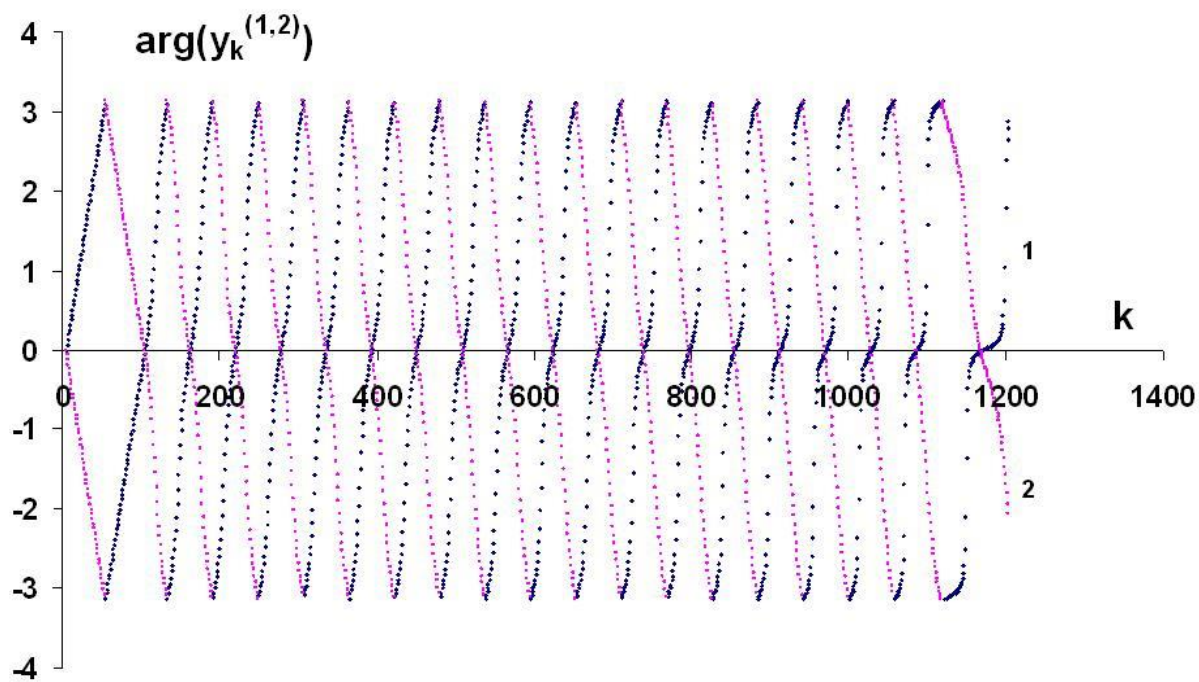

Fig. 5 Graphs of phase of the grid functions $y_k^{(1)}, y_k^{(2)}$ for the wave diffraction on the homogeneous dielectric slab. Matrix $T_k$ is a diagonal.

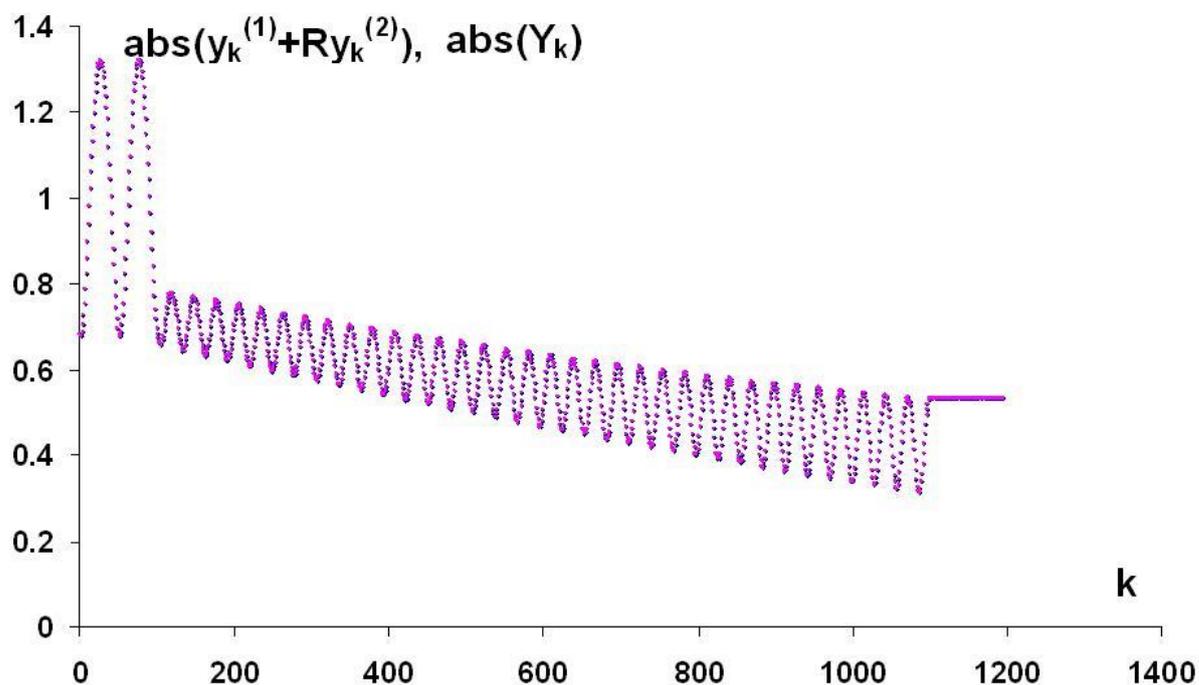

Fig. 6 Graphs of modulus of the grid functions $y_k^{(1)} + R y_k^{(2)}$ and $Y_k$ for the wave diffraction on the homogeneous dielectric slab.

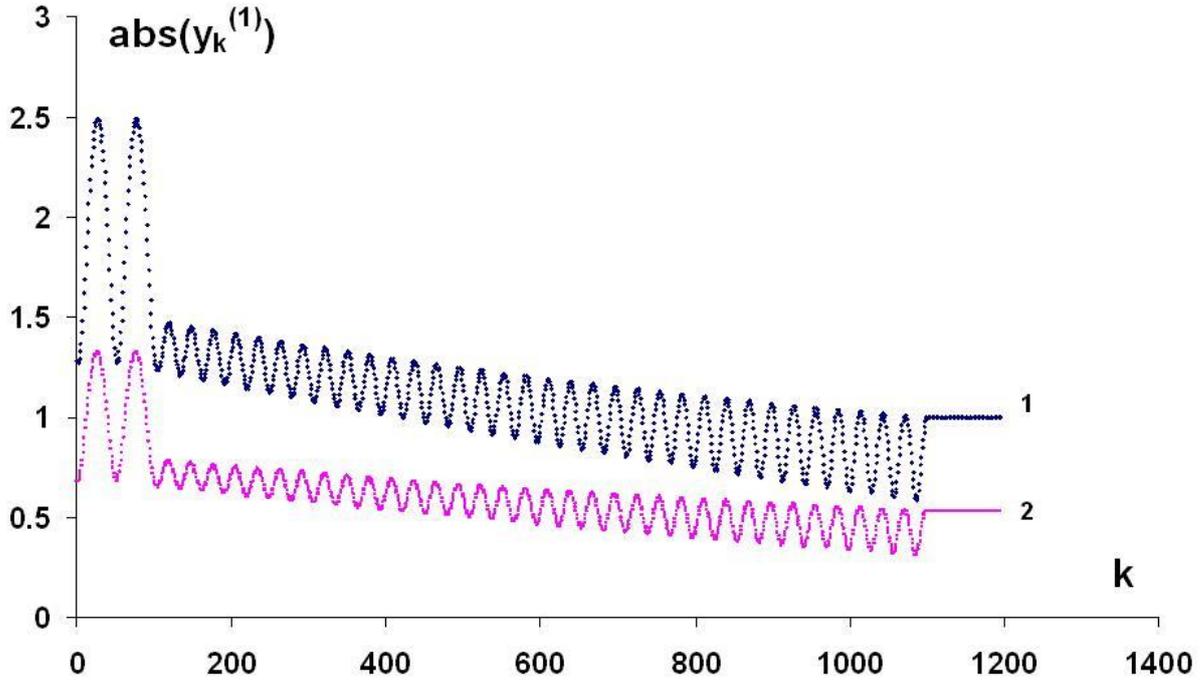

Fig. 7 Graphs of modulus of the grid function $y_k^{(1)}$ for the inverse problem: 1- $y_{N_3}^{(1)} = 1$, 2- $y_{N_3}^{(1)} = 0.531 \times \exp(i \times 1.9865)$

Up to now we have deal with the grid functions $y_k^{(1)}, y_k^{(2)}$ which fulfilled the initial conditions $y_1^{(1)} = 1$, $y_1^{(2)} = 1$ and were calculated on the basis of equations (1.20) with $\rho_k^{(1)}$, $\rho_k^{(2)}$ which are the solutions of the Riccaty equations (1.18) with the initial conditions $\rho_1^{(1)} = \rho_{vac}^{(1)}$, $\rho_1^{(2)} = \rho_{vac}^{(2)}$. In this case before the slab ($1 < k \le 100$) $y_k^{(1)}$ and $y_k^{(2)}$ are the forward and backward waves. But we can consider the inverse case, when $y_k^{(1)}$ and $y_k^{(2)}$ are the forward and backward waves after the slab. If we set the initial conditions $y_{N_3}^{(1)} = 1$, $y_{N_3}^{(2)} = 1$ and $\rho_{N_3}^{(1)} = \rho_{vac}^{(1)}$, $\rho_{N_3}^{(2)} = \rho_{vac}^{(2)}$, and solve the equations (1.20) and the Riccaty equations (1.18) in the inverse order

$$y_k^{(1)} = y_{k+1}^{(1)} / \rho_k^{(1)}$$
$$y_k^{(2)} = y_{k+1}^{(2)} / \rho_k^{(2)} \, . \tag{1.44}$$

$$\rho_k^{(1,2)} = -\frac{b_{k+1}}{\left(\rho_{k+1}^{(1,2)} + a_{k+1}\right)}, \tag{1.45}$$

we obtain the new grid functions that are the independent solutions of the wave equation (1.34), too. After the slab they are the forward and backward waves. So, $y_k^{(2)}$ has no physical meaning and we have to deal with $y_k^{(1)}$ only. Graph of modulus of the grid function $y_k^{(1)}$ for the inverse problem are presented in Fig. 7 (1). From physical point of view, before the slab $y_k^{(1)}$ must be the superposition of the forward and backward waves. Using simple decomposition, we can find amplitudes of the forward and backward waves and find the reflection and transmission coefficients $R^{(inv)} = $ -0.3089-i×9.0980E-002 ($\left| R^{(inv)} \right| = 0.3219$), $T^{(inv)} = $ -0.2445+i×0.4718 ($\left| T^{(inv)} \right| = 0.5314$). Comparison $R^{(inv)}$ and $T^{(inv)}$ with the exact values $R$ and $T$ (or with $R^{(Y)}, T^{(Y)}$) shows good agreement. Setting the initial value $y_{N_3}^{(1)} = Y_{N_3} = 0.531 \times \exp(i \times 1.9865)$ gives the full coincidence of $y_k^{(1,inv)}$ and $Y_k$ (see Fig. 7 (2).).

Results presented above show that the "classic solution" of diffraction problem is one of the two independent solutions of the wave equation that fulfill a certain initial condition. On the bases of the proposed inverse scheme, we can easily find the characteristics of a diffraction problem. Using this scheme instead of S-matrix approach also gives a large gain in computer resources.

### 2.2 Coupling Cavity Model of arbitrary chain of resonators

In the frame of the Coupling Cavity Model (CCM) [13-21] electromagnetic field in each cavity of the chain of resonators are represented as the expansion with the short-circuit resonant cavity modes [41,42,43,44]

$$\vec{E}^{(k)} = \sum_q e_q^{(k)} \vec{E}_q^{(k)}(\vec{r}) \ , \tag{1.46}$$

where $q = \{0, m, n\}$ and such coupling equations for $e_{010}^{(n)}$ can be obtained [13,16,17,20]

$$Z_k e_{010}^{(k)} = \sum_{j=-\infty, j \neq k}^{\infty} e_{010}^{(j)} \alpha_{010}^{(k,j)} . \tag{1.47}$$

Here $e_{010}^{(k)}$ - amplitudes of $E_{010}$ modes, $Z_k = 1 - \dfrac{\omega^2}{\omega_{010}^{(k)2}} - \alpha_{010}^{(k,k)}$, $\omega_{010}^{(k)}$ - eigen frequencies of these modes, $\alpha_{010}^{(k,j)}$ - real coefficients that depend on both the frequency $\omega$ and geometrical sizes of all volumes. Sums in the right side can be truncated

$$Z_k^{(N)} e_{010}^{(N,k)} = \sum_{j=k-N, j \neq k}^{k+N} e_{010}^{(N,j)} \alpha_{010}^{(k,j)} . \tag{1.48}$$

In the case of $N = 1$, the system of coupled equations (1.48) is very similar to the one that can be constructed on the basis of equivalent circuits approach (see, for example [45,46]). But in the frame of the CCM the coefficients $\alpha_{0mn}^{(k,j)}$ are electrodynamically strictly defined for arbitrary $N$ and can be calculated with necessary accuracy. Amplitudes of other modes ( $(m,n) \neq (1,0)$ ) can be found by summing the relevant series

$$e_{0mn}^{(k)} = \frac{\omega_{0mn}^{(k)2}}{\omega_{0mn}^{(k)2} - \omega^2} \sum_{j=k-N}^{k+N} e_{010}^{(j)} \alpha_{0mn}^{(k,j)} . \tag{1.49}$$

Values of the electromagnetic field in the middle of the cavities are:

$$\vec{E}^{(k,m)} \equiv \vec{E}^{(k)}(r=0, z=z_k + d_k/2) = \chi^{(k)} e_{010}^{(k)}, \tag{1.50}$$

where

$$\chi^{(k)} = \left( 1 + \sum_{j=k-N}^{k+N} \frac{e_{010}^{(j)}}{e_{010}^{(k)}} \sum_{\substack{m,n=0,2,4,\dots \\ (m,n)\neq(1,0)}} \frac{\omega_{0mn}^{(k)2} \alpha_{0mn}^{(k,j)}}{\omega_{0mn}^{(k)2} - \omega^2} \right) . \tag{1.51}$$

If we can ignore "long coupling" interaction, the set of coupling equations (1.48) takes the form[6]

$$Z_k e_{010}^{(k)} = e_{010}^{(k-1)} \alpha_{010}^{(k,k-1)} + e_{010}^{(k+1)} \alpha_{010}^{(k,k+1)}, \tag{1.52}$$

where $Z_k = \left( 1 - \dfrac{\omega^2}{\omega_{010}^{(k)2}} - \alpha_{010}^{(k,k)} - i \dfrac{\omega}{\omega_{010}^{(k)} Q_k} \right)$

---

[6] There is a problem of taking into account absorption of RF energy in walls as there are difficulties in obtaining appropriate eigen functions for cylindrical regions. All developed procedures in the frame of the CCM do not include this phenomenon. We used the simplest approach for including absorption into consideration. We supposed that the coupling coefficient do not depend on absorption and include the quality factor into the resonant term in the equations for $e_{010}$ amplitudes.

The equation (1.52) is the second-order difference equation (1.1) ("nondiagonal tight-binding description" [11]) with

$$a_k = -\frac{Z_k}{\alpha_{010}^{(k,k+1)}}$$
$$b_k = \frac{\alpha_{010}^{(k,k-1)}}{\alpha_{010}^{(k,k+1)}} \quad . \qquad . \qquad (1.53)$$
$$f_k = 0$$

Let's find the condition when the transition matrix is the diagonal one. We will consider the $y_k^{(1)}$ as the "forward field", and $y_k^{(2)}$ as the "backward field"[7]. Consider the case when $\rho_k^{(1)} = \rho_* = \exp(i\varphi)$ (the "forward field" has a constant amplitude and a constant phase shift per cell - constant gradient section [47,48]). From (1.26) we obtain that parameters of resonators have to satisfy such equations

$$Z_k = \exp(i\varphi)\alpha_{010}^{(k,k+1)} + \exp(-i\varphi)\alpha_{010}^{(k,k-1)}, \qquad (1.54)$$

and $\rho_k^{(2)} \neq const$ have to be found from the equation

$$\rho_k^{(2)} = -\frac{\alpha_{010}^{(k,k-1)}}{\alpha_{010}^{(k,k+1)}}\frac{1}{\rho_{k-1}^{(2)}} + \exp(i\varphi) + \frac{\alpha_{010}^{(k,k-1)}}{\alpha_{010}^{(k,k+1)}}\exp(-i\varphi), \qquad (1.55)$$

As follow from this expression, a phase shift of the "backward field" is not equal $(-\varphi)$.

As example, consider the chain of cylindrical resonators that are connected via small circular openings in the thin walls. For such chain, we have analytical expressions for the coupling coefficients [41-44]

$$\alpha_{010}^{(k,k)} = -\alpha\frac{r_k^3 + r_{k+1}^3}{R_k^2 d}$$
$$\alpha_{010}^{(k,k+1)} = \alpha\frac{r_{k+1}^3}{R_k^2 d} \qquad , \qquad (1.56)$$

where $r_k$ - the opening radius between $k-1$ and $k$ resonators, $R_k$ - the radius of $k$ cylindrical resonator, $d$ - the resonator length, $\omega_{010}^{(k)} = c\frac{\lambda_{01}}{R_k}$, $\omega = c\frac{\lambda_{01}}{R_*}$, $J_0(\lambda_{01}) = 0$, $\alpha = \frac{2}{3\pi J_1^2(\lambda_{01})}$. If we introduce new notations

$$g_k = \frac{R_k}{R_*}$$
$$u_k = \frac{\alpha r_k^3}{R_*^2 d} \qquad , \qquad (1.57)$$

then real and imaginary parts of (1.54) can be transformed into nonlinear difference equations ($k \geq 1$)

$$g_{k+1}^4 - g_{k+1}^3\frac{(\cos\varphi - 1)}{Q\sin(\varphi)} - g_{k+1}^2 = g_k^4 + g_k^3\frac{(\cos\varphi - 1)}{Q\sin(\varphi)} - g_k^2, \qquad (1.58)$$

---

[7] These components are often referred to as the "forward wave" and the "backward wave". But the notion of "wave" has a standard definition that can not be applied to the arbitrary chain of resonator. So, we will use names "forward field" and the "backward field" to distinguish two solutions.

$$u_k = \frac{1}{2\left(\cos\varphi - 1\right)}\left[ g_k^2 - g_k^4 + g_k^3 \frac{\left(\cos\varphi - 1\right)}{Q\sin(\varphi)} \right], \tag{1.59}$$

Setting the initial value for $u_1$ (proportional to the dimension of the first opening), $g_1$ (proportional to the radius of the first resonator) can be found by solving the quartic

$$g_1^4 - g_1^3 \frac{\left(\cos\varphi - 1\right)}{Q\sin(\varphi)} - g_1^2 + 2\left(\cos\varphi - 1\right)u_1 = 0, \tag{1.60}$$

We will consider the case $\varphi = 2\pi / 3$ [48] and $Q = 10000$.

Dependences of $g_k$ and $u_k$, that support the constant amplitude and phase shift per cell in "forward field", on the resonator number are presented in Fig. 8, Fig. 9 and Fig. 12.

It is interesting to look into the structure of the "backward field" in the constant gradient section[8]. Expression for $\rho_k^{(2)}$ take the following form

$$\rho_k^{(2)} = -\frac{u_k^3}{u_{k+1}^3}\frac{1}{\rho_{k-1}^{(2)}} + \exp(i\varphi) + \frac{u_k^3}{u_{k+1}^3}\exp(-i\varphi), \tag{1.61}$$

Simulation with using a formula (1.20) gives the distribution of the "backward field" along the chain. Analysis of this distribution (see Fig. 10 and Fig. 11) shows that for the most commonly used[9] parameters of chains a phase shift from cell to cell is close to $(-\varphi)$ (see Fig. 11). At the same time, amplitude distributions of the "backward field" are the growing one and amplitude differences between the start and the end can be significant (see Fig. 10).

For small[10] initial value $u_1$, phase distribution differs substantially from the law $(-k\varphi)$ and amplitudes distribution is strongly nonuniform. It can create some difficulties in the process of tuning such structures.

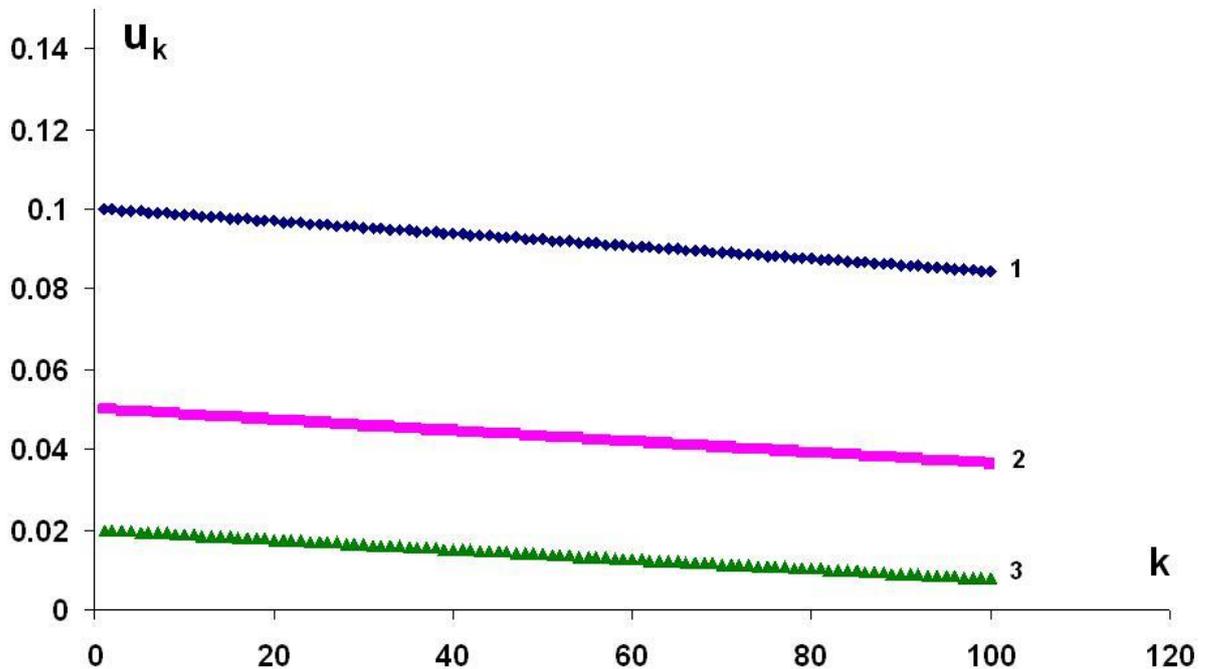

Fig. 8 $u_k$ versus the resonator number for different $u_1$

---

[8] It is important to know structure of this field for developing the tuning methods (see, for example, [15,48])

[9] For the SLAC structure $u_1 \approx 0.03$ [48]

[10] It can be realized in the short structures that are developed for study electrical breakdown in high-gradient accelerating structures

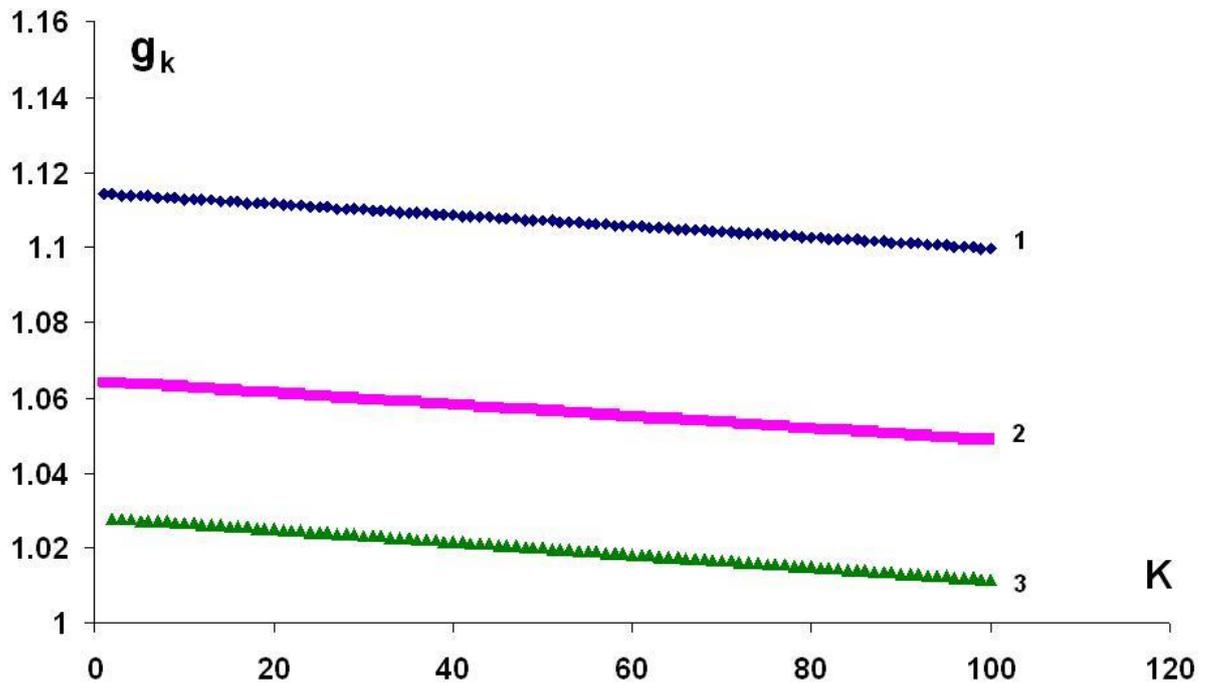

Fig. 9  $g_k$ versus the resonator number for different $u_1$ (1- $u_1$=0.1, 2- $u_1$=0.05, 3- $u_1$=0.02)

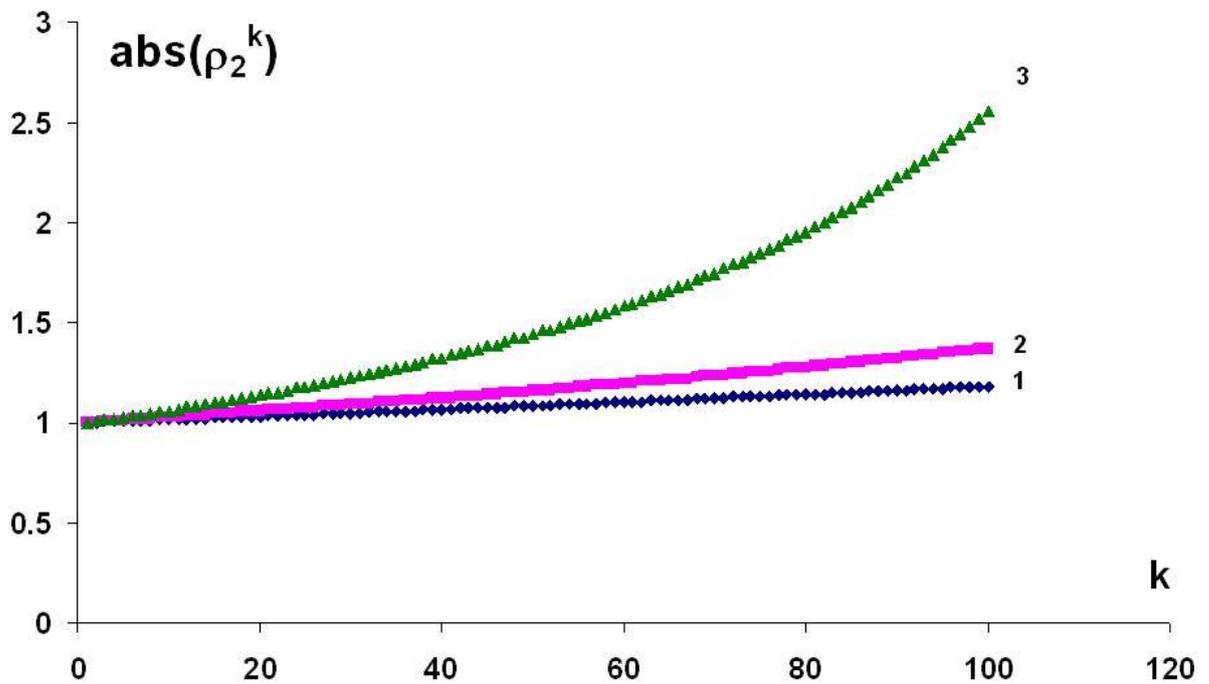

Fig. 10 Amplitude distribution of the "backward field" for different $u_1$ (1- $u_1$=0.1, 2- $u_1$=0.05, 3- $u_1$=0.02)

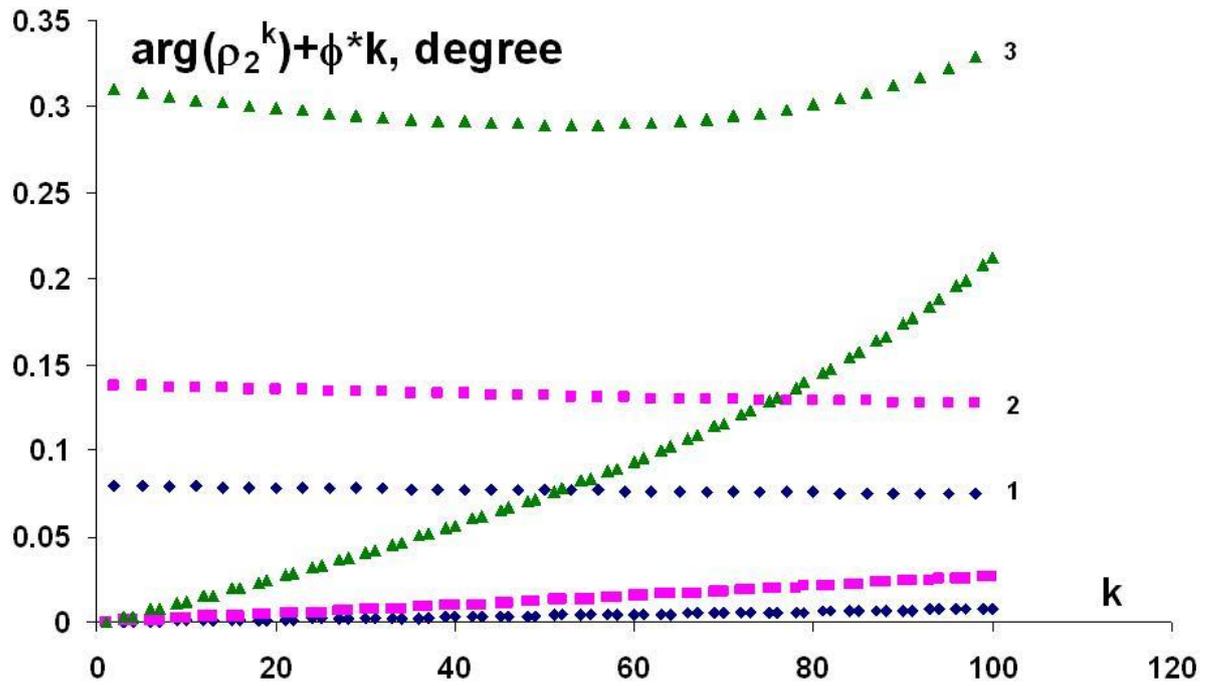

Fig. 11 Phase deviation of the "backward field" from the law $(-k\varphi)$ for different $u_1$ (1- $u_1$ =0.1, 2- $u_1$ =0.05, 3- $u_1$ =0.02)

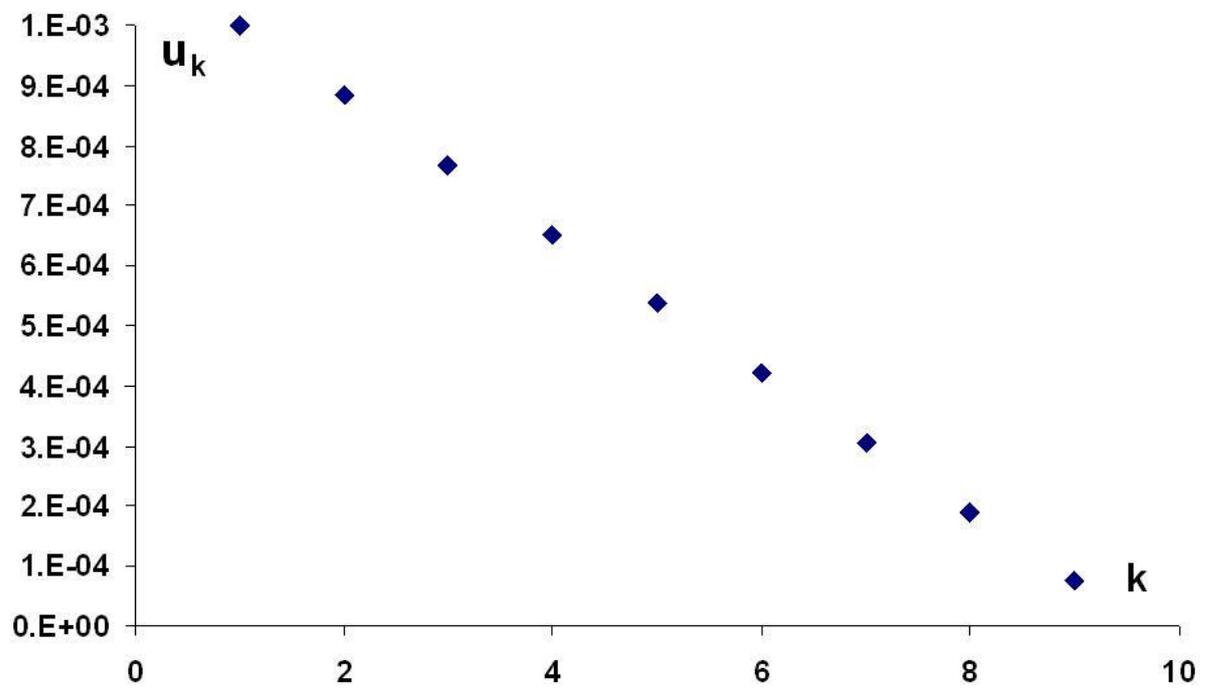

Fig. 12 $u_k$ versus the resonator number for small value of $u_1$

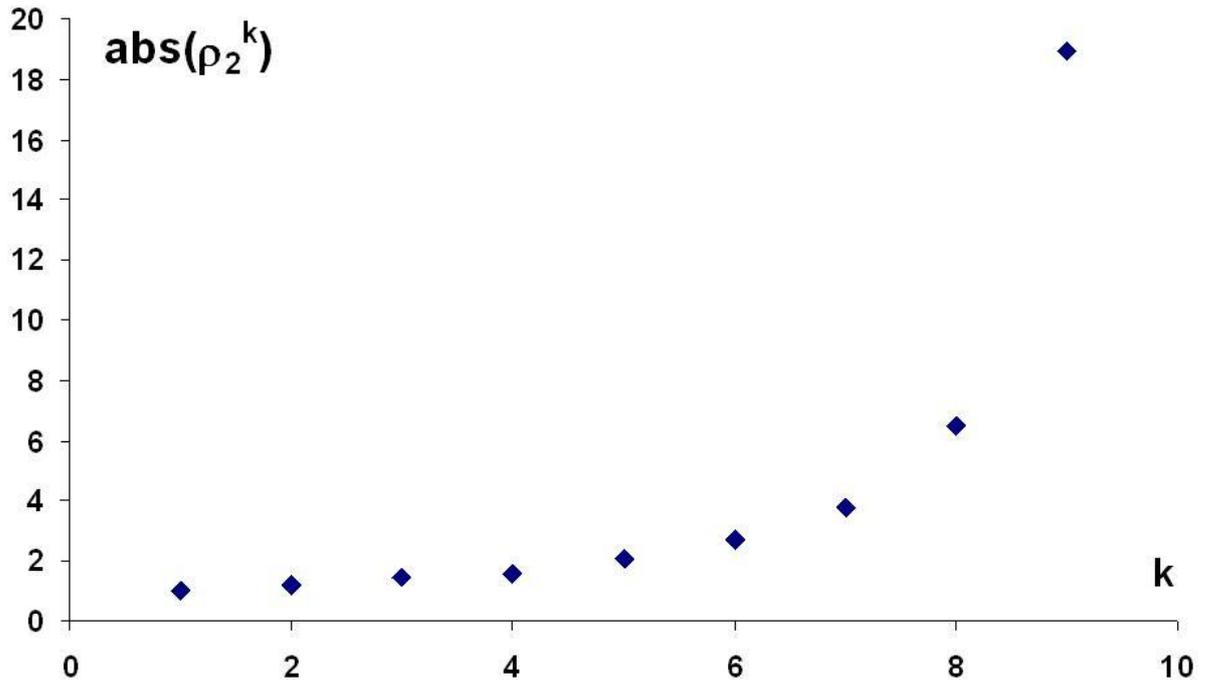

Fig. 13 Amplitude distribution of the "backward field" for small value of $u_1 = 0.001$

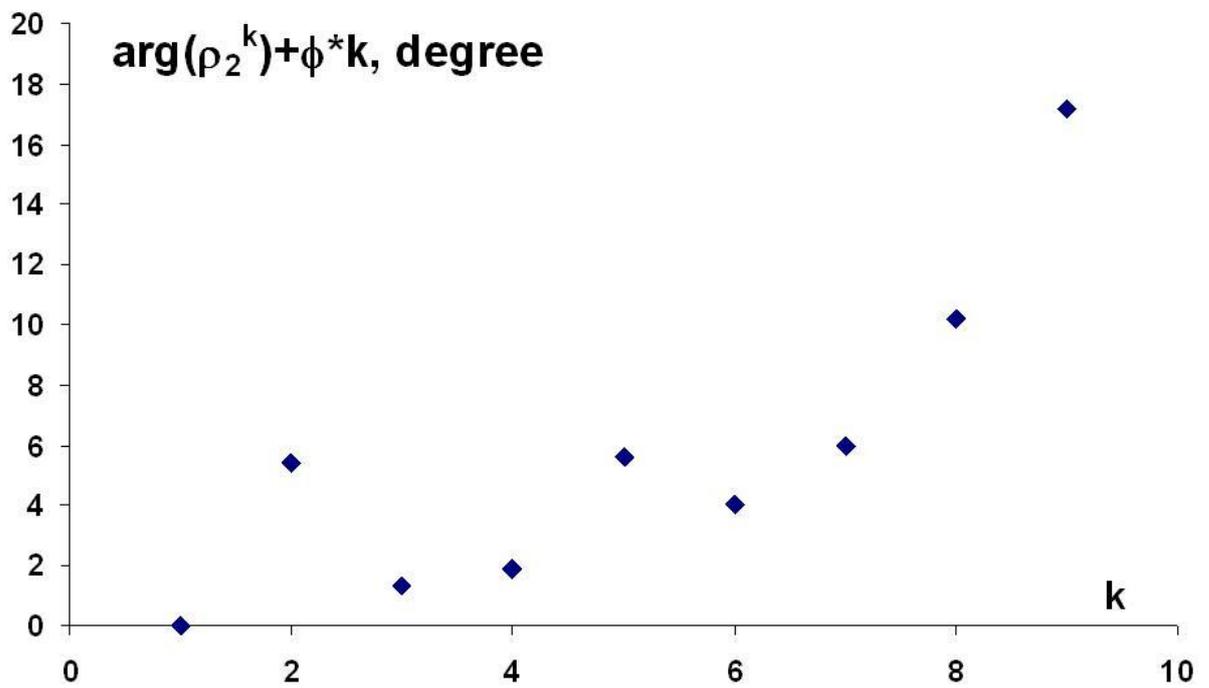

Fig. 14 Phase deviation of the "backward field" from the law $(-k\varphi)$ for small value of $u_1 = 0.001$

## Conclusions

We presented results of our research of possibility of transformation of a difference equation into a system of the first-order difference equations. Several examples show that proposed approach can be useful in solving different physical problems.

## Acknowledgments

I would like to thank M.M.Aizatskyi for useful discussion of the results.

## References

1 Sedaghat Hassan. Form Symmetries and Reduction of Order in Difference Equations. CRC Press, 2011

2 Brekhovskikh L.M. Waves in Layered Media.. Academic Press,  Elsevier, 1960

3 Ginzburg V.L. The Propagation of Electromagnetic Waves in Plasmas. Pergamon, London; Addison-Wesley, Reading, Mass., 1964

4 R. Jacobsson Light reflection from films of continuously varying refractive index. Progress in Optics, Vol. 5, ed. E. Wolf . Amsterdam, North Holland Publishing Company and New York, J. Wiley and Sons, 1965

5 Au Kong Jin. Electromagnetic Wave Theory. John Wiley & Sons

6 Chew  Weng Cho. Waves and Fields in Inhomogenous Media. Wiley-IEEE Press, 1999

7 Born Max, Wolf Emil. Principles of optics -  electromagnetic theory of propagation, interference and diffraction of light. Cambridge University Press, 2005

8 Yeh Pochi. Optical Waves in Layered Media. Wiley, 2005

9 Barclay Les. Propagation of Radiowaves. The Institution of Engineering and Technology,. 2013

10 Ruey-Bing (Raybeam) Hwang Periodic structures : mode-matching approach and applications in electromagnetic engineering, John Wiley & Sons Singapore Pte.Ltd., 2013

11 Akkermans E., Dunne G. V., Levy E.. Wave Propagation in One Dimension: Methods and Applications to Complex and Fractal Structures. In: Negro L.D. (ed.) Optics of aperiodic structures. Fundamentals and Device Applications. Taylor & Francis Group, 2014

12 Banerjee P.P., Jarem J.M. Computational Methods for Electromagnetic and Optical Systems CRC Press, 2014

13 M.I.Ayzatsky. New Mathematical Model of an Infinite Cavity Chain. Proceedings of the EPAC96, 1996,v.3, p.2026-2028; On the Problem of the Coupled Cavity Chain Characteristic Calculations.  http://xxx.lanl.gov/pdf/acc-phys/9603001.pdf, LANL.arXiv.org e-print archives, 1996

14 M.I.Ayzatskiy, K.Kramarenko. Coupling coefficients in inhomogeneous cavity chain Proceedings of the EPAC2004, 2004, pp.2759-2761

15 M.I. Ayzatskiy, V.V. Mytrochenko. Electromagnetic fields in nonuniform disk-loaded waveguides. http://lanl.arxiv.org/ftp/arxiv/papers/1503/1503.05006.pdf, LANL.arXiv.org e-print archives, 2015; M.I. Ayzatskiy, V.V. Mytrochenko. Electromagnetic fields in nonuniform disk-loaded waveguides. PAST, 2016, N.3, pp.3-10

16 M.I. Ayzatskiy, V.V. Mytrochenko. Coupled cavity model based on the mode matching technique. http://lanl.arxiv.org/ftp/arxiv/papers/1505/1505.03223.pdf, LANL.arXiv.org e-print archives, 2015

17 M.I. Ayzatskiy, V.V. Mytrochenko. Coupled cavity model for disc-loaded waveguides. http://lanl.arxiv.org/ftp/arxiv/papers/1511/1511.03093.pdf, LANL.arXiv.org e-print archives, 2015 18 M. I. Ayzatsky, V. V. Mytrochenko. Numerical design of nonuniform disk-loaded waveguides. http://lanl.arxiv.org/ftp/arxiv/papers/1604/ 1604.05511.pdf, LANL.arXiv.org e-print archives, 2016

19 M.I. Ayzatskiy, V.V. Mytrochenko. Numerical investigation of tuning method for nonuniform disk-loaded waveguides. http://lanl.arxiv.org/ftp/arxiv/papers/1606/ 1606.04292.pdf, LANL.arXiv.org e-print archives, 2016

20 M.I. Ayzatskiy, V.V. Mytrochenko. Methods for calculation of the coupling coefficients in the Coupling Cavity Model of arbitrary chain of resonators. http://lanl.arxiv.org/ftp/arxiv/papers//1609/1609.01481.pdf, LANL.arXiv.org e-print archives, 2016

21 M.I. Ayzatskiy, V.V. Mytrochenko. Coupling cavity model for circular cylindrical waveguide with uniform cross section. http://lanl.arxiv.org/ftp/arxiv/papers//1609/ 1609.03306. pdf, LANL.arXiv.org e-print archives, 2016

22 Y. Liu , B. Vial, S. A. R. Horsley et al. Direct manipulation of wave amplitude and phase through inverse design of isotropic media, http://xxx.lanl.gov/pdf/1701.06021v1, LANL.arXiv.org e-print archives, 2017

23 L. Beilina, E. Smolkin Computational design of acoustic materials using an adaptive optimization algorithm http://xxx.lanl.gov/pdf/1701.06006v, LANL.arXiv.org e-print archives, 2017

24 Jordan Charles. Calculus of Finite Differences. Chelsea Publishing Company, 1950

25 Cull  Paul, Flahive Mary, Robson Robby Difference Equations. From Rabbits to Chaos. Springer, 2005

26 Levy H., Lessman F. Finite difference equations. Macmillan, 1961

27 Brand Louis. Differential and Difference Equations. John Wiley & Sons, 1966

28 Miller Kenneth S. An introduction to the calculus of finite differences and difference equations. Dover Publications, 1966

29 Wimp Jet. Computation with recurrence relations. Pitman Publishing, 1984

30 Kocic V.L., Ladas G. Global behavior of nonlinear difference equations of higher order with applications. Kluwer, 1993

31 Kulenovic Mustafa R.S., Ladas G. Dynamics of second order rational difference equations. Chapman and Hall_CRC, 2001

32 Lakshmikantham  V., Trigiante Donato. Theory of Difference Equations. Numerical Methods and Applications. CRC Press, 2002

33 Ashyralyev Allaberen , Sobolevskii Pavel E. New Difference Schemes for Partial Differential Equations.. Springer Basel AG, 2004

34 Grove E.A., Ladas G. Periodicities in nonlinear difference equations. CRC, 2005

35 Sedaghat Hassan. Nonlinear Difference Equations. Theory with Applications to Social Science. Models-Springer Netherlands, 2003

36 Banasiak J. Mathematical Modelling in One Dimension. An Introduction via Difference and Differential Equations. Cambridge University Press, 2013

37 Saber Elaydi. An Introduction to Difference Equations. Springer, 2005

38 Samarskii A.A., Nikolaev E.S. Numerical Methods for Grid Equations. Volume I Direct Methods. Birkhäuser Basel (1989)

39 Dobrowolski Janusz A. Microwave Network Design Using the Scattering Matrix. ARTECH HOUSE, 2010

40 Jackson I.D. Classical_electrodynamics. John Wiley&Sons, 1962

41 Bethe H.A. Bathe H.A. Theory of Diffraction by Small Holes. Phys. Rev., 1944, v.66, N7, p.163-182.

42 V.V. Vladimirsky. ZhTF, 1947, v.17, N.11, p.1277-1282. Владимирский В.В. Связь полых электромагнитных резонаторов через малое отверстие. ЖТФ, 1947, т.17, №11, с.1277-1282 43 A.I. Akhiezer, Ya.B. Fainberg. UFN, 1951, v.44, N.3, p.321-368.Ахиезер А.И., Файнберг Я.Б. Медленные волны. УФН, 1951, т.44, №3, с.321-368

44 R.M. Bevensee. Electromagnetic Slow Wave Systems. John Wiley\&Sons, Inc.,New York-London-Sydney, 1964

45 D.E. Nagle, E.A. Knapp, B. C. Knapp, Coupled Resonator Model for Standing Wave Accelerator Tanks. The Review of Scientific Instruments, 1967, V.38, N.11, pp.1583-1587

46 D.H.Whittum Introduction to Electrodynamics for Microwave Linear Accelerators. In: S.I.Kurokawa, M.Month, S.Turner (Eds) Frontiers of Accelerator Technology, World Scientific Publishing Co.Pte.Ltd., 1999

47 R.B.Neal Theory of the constant gradient linear electron accelerator. M.I.Report No.513, 1958

48 R.B Neal, General Editor, The Stanford Two-Mile Accelerator, New York, W.A. Benjamin, 1968